\documentclass[12pt]{amsart}

\usepackage{amscd,amssymb}
\usepackage[all]{xy}
\usepackage[plainpages,backref,urlcolor=blue]{hyperref}
\usepackage[utf8]{inputenc}
\usepackage{mathtools}

\theoremstyle{plain}
\newtheorem{thm}[subsection]{Theorem}
\newtheorem{lem}[subsection]{Lemma}
\newtheorem{prop}[subsection]{Proposition}

\theoremstyle{definition}
\newtheorem{rk}[subsection]{Remark}

\newtheorem{ex}[subsection]{Example}

\numberwithin{equation}{section}
\newcommand{\al}{\alpha}
\newcommand{\be}{\beta}

\newcommand{\C}{\mathbb C}

\DeclareMathOperator{\Ann}{Ann}
\DeclareMathOperator{\Fitt}{Fitt}
\DeclareMathOperator{\Syz}{Syz}
\DeclareMathOperator{\ord}{ord}

\begin{document}

\title[Plane curve singularities and Fitting ideals]{Plane curve singularities and Fitting ideals}

\author[Alexandru Dimca]{Alexandru Dimca$^1$}
\address{Universit\'e C\^ote d'Azur, CNRS, LJAD, France and Simion Stoilow Institute of Mathematics,
P.O. Box 1-764, RO-014700 Bucharest, Romania}
\email{Alexandru.Dimca@univ-cotedazur.fr}

\author[Gabriel Sticlaru]{Gabriel Sticlaru}
\address{Faculty of Mathematics and Informatics,
Ovidius University
Bd. Mamaia 124, 900527 Constanta,
Romania}
\email{gabriel.sticlaru@gmail.com }

\thanks{$^1$ partial support from the project ``Singularities and Applications'' - CF 132/31.07.2023 funded by the European Union - NextGenerationEU - through Romania's National Recovery and Resilience Plan.}



\begin{abstract}
In this note we investigate  the Fitting ideals associated to the Tjurina ideal of a non quasi-homogeneous plane curve singularity.
Special properties occur when the difference between Milnor number and Tjurina number is at most 2.

\end{abstract}

\maketitle

\section{Introduction}

Let $R=\C\{x,y\}$ be the convergent power series ring in $2$ variables $x$ and $y$ with complex coefficients, and let $X:f=0$ be an isolated curve singularity at the origin $0 \in \C^2$, with $f \in R$.
We denote by $J_f$ the Jacobian ideal of $f$, i.e. the ideal in $R$ spanned by the partial derivatives 
$$f_x=\partial_{x}f \text{ and } f_y=\partial_{y}f.$$
 and  by $M(f)=R/J_f$ the corresponding  quotient ring, called the  Milnor algebra of $f$.
  We denote by $I_f$ the Tjurina ideal of $f$, i.e. the ideal in $R$ spanned by 
$f_x=\partial_{x}f , \ f_y=\partial_{y}f$ and $f$,
 and  by $T(f)=R/I_f$ the corresponding  quotient ring, called the Tjurina algebra of $f$.
Then both $M(f)$ and $T(f)$ are Artinian local $\C$-algebras, and
$$\mu(f)= \dim_{\C} M(f) \text{ and } \tau(f)= \dim_{\C} T(f)$$
are the Milnor and resp. Tjurina number of $f$.

It is well known that $X:f=0$ is a quasi-homogeneous singularity if and only if $I_f=J_f$, or equivalently $\mu(f)=\tau(f)$, see \cite{S}.
In this case the minimal resolution of $I_f$, which is a complete intersection, has the following form
\begin{equation}
\label{res0}
 0 \to R  \stackrel{u'} \longrightarrow R^2  \stackrel{v'} \longrightarrow I_f \to 0,
 \end{equation} 
where $u'(1)=(-f_y,f_x)$ and $v'(1,0)=f_x$, $v'(0,1)=f_y$.

{\it We assume from now on in this note that $X:f=0$ is not a quasi-homogeneous singularity. }
Let
$$\Syz(f)= \ker v= \{(a,b,c) \in R^3 \ : \ af_x+bf_y+cf=0\}$$
be the first syzygy module of the $R$-module $I_f$, and 
$$D(f)=\{ \delta=a\partial_x + b \partial_y \in D(R) \  : \ \delta (f) \in  Rf \},$$
where $D(R)$ is the $R$ module of $C$-derivations of $R$, and $Rf$ is the principal ideal generated by $f$ in $R$.
Then $D(f)$ is the $R$-module of logarithmic vector fields tangent to $f$ and there is an obvious isomorphism
$$\Syz(f) = D(f) \text{ given by } (a,b,c) \mapsto a\partial_x + b \partial_y.$$
Moreover $D(R)$ is a free $R$-module of rank 2, see \cite[Corollary 1.7]{S}.Then one has the following result, in which the first three claims restate our discussion above and recall the definition of the Fitting ideals.

\begin{prop}
\label{prop1}
With the above notation,  one has the following.
\begin{enumerate}

\item The $R$-module $\Syz(f)$ is free of rank 2.

\item If $\rho_1$ and $\rho_2$ are a basis of the $R$-module $\Syz(f)$,
then 
$$0 \to R^2 \stackrel{u} \longrightarrow R^3 \stackrel{v} \longrightarrow I_f \to 0$$
is a minimal free resolution for $I_f$, where
$$u(s,t)=s\rho_1+t\rho_2$$
 and
 $$v(a,b,c)=af_x+bf_y+cf.$$

\item Let  $U$ denote the matrix associated to the morphism $u$ and let
$I_k(U)$ be the ideal in $R$ generated by all the $k$-minors in $U$. 
Then  $I_j(U)=\Fitt_{3-j}(I_f)$ for $j=1,2$, the corresponding Fitting ideals of the $R$-module $I_f$.

\item In addition, $I_f=I_2(U)=\Fitt_1(I_f)$ and in particular
$$\tau(f)=\dim_{\C} \frac{R}{\Fitt_1(I_f)}.$$

\item The ideal
$$\Ann(f)= \{g \in R \ : \ gf=0 \text{ in } M(f) \} \subset R$$
is minimally generated by exactly two elements, and hence it is a 0-dimensional complete intersection.

\end{enumerate}

\end{prop}

The next result tells us something about the other Fitting ideal $\Fitt_2(I_f)=I_1(U)$, generated by all the entries of the matrix $U$.

\begin{prop}
\label{prop2}
With the above notation, one has the following.
\begin{enumerate}

\item $$\Ann(f) \subset \Fitt_2(I_f) \subset (x,y)={\bf m}.$$

\item $$\mu(f)-\tau(f) \geq \dim_{\C} \frac{R}{\Fitt_2(I_f)}.$$

\end{enumerate}
In particular, the equalities hold in (1) and (2) if $\mu(f)-\tau(f)=1$.

\end{prop}

 The  case  $\mu(f)-\tau(f)=2$ is settled in the following result.

\begin{thm}
\label{thm3}
With the above notation, if $\mu(f)-\tau(f)=2$, then
$$\Ann(f) = \Fitt_2(I_f) \text{  and }
\mu(f)-\tau(f)= \dim_{\C} \frac{R}{\Fitt_2(I_f)}.$$

\end{thm}
 Example \ref{ex4}  shows that the equalities in Theorem \ref{thm3} may fail when the codimension of $\Ann(f)$ is at least 3.

\section{Proof of Proposition \ref{prop1}}

The claims (1), (2) and (3) are clear from the discussion proceeding this result.
Using now Hilbert-Burch Theorem, see \cite[Theorem 20.15]{Eis0}, it follows that
$I_f=aI_2(U)$, where $a \in R$. Since $I_f$ has finite codimension in $R$, it follows that $a$ has to be a unit, and hence we get the claim (4).

To prove (5), let
\begin{equation}
\label{syz1}
\rho_j=(a_j,b_j,c_j)
\end{equation}
for $j=1,2$ be a basis of the free module $\Syz(f)$. Then $g \in \Ann(f)$ if and only if there exists $a',b' \in R$ such that
$$a'f_x+b'f_y+gf=0.$$
This means that $(a',b',g) \in \Syz(f)$, and hence 
this happens if and only if $g \in (c_1,c_2)$, the ideal generated by $c_1$ and $c_2$ in $R$.
Suppose that $\Ann(f)=(c_1,c_2)$ is a principal ideal in $R$, say generated by $c \in R$. Then $c_1$ and $c_2$ are divisible by $c$, and this shows that a basis of the free module $\Syz$ can have the following form
$$\rho_1'=(a,b,c) \text{ and } \rho_2'=(-f_y,f_x,0).$$
Then the ideal $\Fitt_1(I_f)=I_2(U)$ is generated by $cf_x$,  $cf_y$ and
$af_x+bf_y=-cf$. This implies that $I_f=I_2(U) \subset (c)$,
which is a contradiction, since the principal ideal $(c)$ has infinite codimension 
(as a vector subspace) in $R$.

\begin{rk}
\label{rk1}
If the basis of the free $R$-module $\Syz(f)$ is chosen as in \eqref{syz1}, it follows from Hilbert-Burch Theorem, see \cite[Theorem 20.15]{Eis0}, that there is a unit $h \in R$ such that
$$a_1b_2-a_2b_1=hf, \ a_1c_2-a_2c_1=-hf_y \text{ and } b_1c_2-b_2c_1=hf_x.$$
\end{rk}
\begin{rk}
\label{rk2}
There seems to be no similar result for germs of isolated surface singularities in $\C^3$. For instance, if we take 
$$f=x^6+y^5+z^5+x^3y^2z+z^7 \in R'=\C\{x,y,z\}$$
a direct computation with SINGULAR shows that
$\mu(f)=80$, $\tau(f)=68$ and the module $\Syz(f)$ has 6 generators.
If we denote by $U$ the corresponding $4 \times 6$ matrix, then
$$\dim_{\C} \frac{R'}{I_1(U)}=11, \ \dim_{\C} \frac{R'}{I_2(U)}=89, 
\ \dim_{\C} \frac{R'}{I_3(U)}=312 \text{ and } I_4(U)=0.$$
On the other hand $\Ann(f)=(x^2,5y^2-yz^2,z^3)$ and of course
$$ \dim_{\C} \frac{R'}{\Ann(f)}=\mu(f)-\tau(f)=12,$$
since the exact sequence \eqref{esAnn} and the equality \eqref{eqAnn} below clearly hold in any dimension.
In this example $\Ann(f)$ is still a 0-dimensional complete intersection, but this property fails in general. Indeed, for the singularity
$$f=x^5+xy^5+yz^6+x^2y^3z^3+x^3z^5$$
a direct computation with SINGULAR shows that
$$\Ann(f)=(6y^2+y^3z-xz^3,
x^3, 10x^2y+2z^4+x^2y^2z,
120xz^2+239y^4z, yz^2, z^6).$$

\end{rk}

\begin{ex}
\label{ex1}
Let $f=x^2y^2+x^5+y^5$. Then a direct computation using SINGULAR \cite{Sing} shows that
$$\rho_1=(8xy-75x^2y^2,-8x^2+10y^3+50x^3y,20x-125x^2y)$$
and
$$\rho_2=(12y^2+10x^3-75xy^3,-12xy+50x^2y^2,20y-125xy^2)$$
give a basis for the free module $\Syz(f)$.
It follows in particular that
$$\Ann(f)=(c_1,c_2)=(20x-125x2y,20y-125xy^2)=(x,y)={\bf m},$$
the maximal ideal in $R$.
It is known that $\mu(f)= 11$ and $\tau(f)= 10$, see \cite[Example (6.56)]{RCS}.
\end{ex}

\section{Proof of Proposition \ref{prop2}}

The inclusion $\Ann(f) \subset I_1(U)$ follows from our remark in the proof of Proposition \ref{prop1} that $\Ann(f)=(c_1,c_2)$. 
The inclusion $I_1(U) \subset {\bf m}$ can be proved as follows.
First one has $c_1 \in {\bf m}$ and $c_2 \in {\bf m}$, since otherwise $X$ would be a quasi-homogeneous singularity by \cite{S}.
On the other hand, the coefficients $a_1,b_1,a_2$ and $b_2$ are also in
${\bf m}$. Indeed, suppose that $a_1$ is an invertible element in $R$.
Then $f_x \in (f_y,f)$, which is a contradiction. Indeed, if we write
\begin{equation}
\label{syz2}
f=A_p(y)x^p + \text{ higher order terms in } x \in R=\C\{y\}\{x\},
\end{equation}
with $A_p(y) \ne 0$, then $f_x$ has order $p-1$ in $x$, but $f$ and $f_y$ have both order $\geq p$ in $x$.

To prove (2) we use first the exact sequence
\begin{equation}
\label{esAnn}
0 \to \frac{\Ann(f)}{J_f} \to M(f) \stackrel{f} \longrightarrow M(f) \to T(f) \to 0
\end{equation}
to get
$$\dim_{\C} \frac{\Ann(f)}{J_f}=\tau(f).$$
It follows that
\begin{equation}
\label{eqAnn}
\tau(f)=\mu(f) - \dim_{\C} \frac{R}{\Ann(f)}
\end{equation}
and hence
$$
\tau(f)  \leq \mu(f)-\dim_{\C} \frac{R}{\Fitt_2(I_f)},$$
where the last inequality follows from (1).

When $\mu(f)-\tau(f)=1$, it follows that the class of $f$ in $I_f/J_f$ yields a basis of this 1-dimensional vector space. Hence
$\Ann(f) = \bf m$, and this gives the required equalities in (1) and (2).

\endproof

\begin{ex}
\label{ex2}
In general, the difference 
$$\Delta(f)=\mu(f)-\tau(f) -\dim_{\C} \frac{R}{\Fitt_2(I_f)}\geq 0$$
can be quite large. For instance, if we take
$$f= 2y^{14}+x^{15}+2x^{13}y^2+x^{11}y^4,$$
we get by a direct computation with SINGULAR \cite{Sing}
$$
\mu(f) = 182, \  \tau(f) = 162 \text{ and } \dim_{\C} \frac{R}{\Fitt_2(I_f)}=6.$$
It follows that $\Delta(f)=14$ in this case.
\end{ex}

\section{Proof of Theorem \ref{thm3}}

We have seen in Proposition \ref{prop1} (5) that the ideal $\Ann(f)$ is a complete intersection at the origin $0 \in \C^2$. From the classification of the simplest such complete intersections, see for instance \cite[Chapter 9]{RCS}, we have the following basic result.
\begin{lem}
\label{lem1}
With the above notation, the following choices for $c_1$ and $c_2$ in 
\eqref{syz1} are the only  ones, up to a change of basis for $\Syz$ and an analytic change of coordinates at the origin of $\C^2$.
\begin{enumerate}

\item If  $\dim_{\C} \frac{R}{ \Ann(f) }=1$, then $c_1=x$ and $c_2=y$.

\item If  $\dim_{\C} \frac{R}{ \Ann(f) }=2$, then $c_1=x$ and $c_2=y^2$.

\item If  $\dim_{\C} \frac{R}{ \Ann(f) }=3$, then $c_1=x$ and $c_2=y^3$.

\end{enumerate}

\end{lem}

Note that for  $\dim_{\C} \frac{R}{ \Ann(f) }=4$, the pair $(c_1,c_2)$ is no longer determined. Indeed, one may have  $(c_1,c_2)=(x,y^4)$ or
$(c_1,c_2)=(x^2,y^2)$.

Now we start the proof of Theorem \ref{thm3} assuming that
$(c_1,c_2)=(x,y^2)$, which is possible by Lemma \ref{lem1}. It is enough to prove the equality $\Ann(f) = I_1(U)$. Note that one clearly has by Proposition \ref{prop2} (1) and Lemma \ref{lem1} (2) the inclusions
$$(x,y^2)=(c_1,c_2)=\Ann(f) \subset I_1(U) \subset (x,y).$$
To prove that $\Ann(f) = I_1(U)$ is equivalent to show that $I_1(U) \ne (x,y)$, since there is no ideal $K$ such that $(x,y^2) \subsetneq K \subsetneq (x,y)$. Since $c_1=x$ and $c_2=y^2$, the equalities in Remark \ref{rk1} yield the following
\begin{equation}
\label{syz10}
y^2a_1-xa_2=-hf_y \text{ and } y^2b_1-xb_2=hf_x
\end{equation}
for some unit $h \in R$. 
Note also that $I_1(U) = (x,y)$ if and only if $y$ occurs with a non-zero coefficient in at least one of the series $a_1, a_2,b_1$ and $b_2$.

\medskip

\noindent {\bf Case 1} $\ord f \leq 2$.
If $\ord f \leq 2$, then $X$ is either smooth, or a singularity of type $A_k$ for some $k$, see \cite{AGV, RCS} for the classification of simple singularities $A_k$, $D_k$, $E_6$, $E_7$ and $E_8$. Since all these simple singularities are quasi-homogeneous,  the singularity $X$ is quasi-homogeneous in this case. It follows that this case cannot occur.

\medskip

\noindent {\bf Case 2} $\ord f = 3$.
 If we write
$$f=f_3+f_4+ f_5 + \ldots$$
for the decomposition of $f$ into homogeneous components, we see that
$f$ is a simple singularity of type $D_k$, unless $f_3 \ne 0$ is a perfect cube, let's say $f_3=\ell^3$ for some linear form $\ell$. To prove this claim, one may have a look at \cite[Proposition 8.17]{RCS}.

Then \eqref{syz10} imply that $\ell$ as well as the linear parts of $a_2$ and $b_2$ are divisible by $x$. Hence 
$$f_3= \al x^3 \text{ for some } \al \in \C^*.$$ 
Suppose that $I_1(U) = (x,y)$. Then one of the following two situations has to occur.

\medskip

\noindent {\bf Subcase 2A} $a_1=\al' y + \be' x+ r_2$, with $r_2 \in (x,y)^2$ and $\al' \ne 0$.
The first relation in \eqref{syz10} implies that $y^3$ occurs in the product $-hf_y$. It follows that $y^3$ occur in $f_y$ (resp. $y^4$ occurs in $f$) with a non-zero coefficient. In other words
$$f= \al x^3+ \be y^4 + xg+h,$$
where $g \in (x,y)^3$ and $h \in (x,y)^5$. Then it follows that
$f$ has type $E_6$, see for instance the proof of \cite[Proposition 8.21]{RCS}.

\medskip

\noindent {\bf Subcase 2B}  Assume that we are not in Subase 2A, namely the linear part of $a_1$ does not contain $y$. Then necessarily $b_1=\al' y + \be' x + r_2$, with $r_2\in (x,y)^2$ and $\al' \ne 0$, $\be \in \C$.

The second relation in \eqref{syz10} implies that $y^3$ occurs in the product $hf_x$. It follows that  $y^3$ occur in $f_x$, and hence $xy^3$ occurs in $f$ with a non-zero coefficient. If we write
$$f=\al x^3 + \be xy^3+ x^2g+h$$
where $g \in (x,y)^2$ and $h \in (x,y)^5$, we see that $f$ is a simple singularity of type $E_7$, as in the proof of \cite[Proposition 8.21]{RCS}.

In all this cases $X:f=0$ is a quasi-homogeneous singularity, and hence
$\Ann(f)=R$, a contradiction. Hence our claim is proved for $\ord f \leq 3$.

\medskip

\noindent {\bf Case 3} $\ord f =4$.
Now we have
$$f=f_4+ f_5 + \ldots$$
with $f_4 \ne 0$ for the decomposition of $f$ into homogeneous components. Note that \eqref{syz10} implies that the linear parts 
$(a_2)_1$ of $a_2$ and $(b_2)_1$ of $b_2$ vanish. Hence if we take the homogeneous component of degree 3 in the first equality from Remark \ref{rk1} we get the following equality
\begin{equation}
\label{eqK1}
(a_1)_1(b_2)_2-(a_2)_2(b_1)_1=0,
\end{equation}
where $g_2$ means the homogeneous part of degree 2 in some $g \in R$.
Next we take the homogeneous components of degree 3 in the equalities from \eqref{syz10} and  get the following 
\begin{equation}
\label{eqK2}
y^2(a_1)_1-x(a_2)_2=-\al(f_4)_y \text{ and } y^2(b_1)_1-x(b_2)_2=\al(f_4)_x ,
\end{equation}
where $\al=h(0) \ne 0$. Multiply now the first equality by $-(b_1)_1$, the second by $(a_1)_1$ and make the sum. In view of \eqref{eqK2} we get
\begin{equation}
\label{eqK3}
(a_1)_1(f_4)_x+(b_1)_1(f_4)_y=0.
\end{equation}
We assume from now on in this Case that $I_1(U)=(x,y)$ and get a contradiction. With this assumption \eqref{eqK3} is a non-trivial Jacobian syzygy for $f_4 \ne 0$. Let $r$ be the number of distinct factors in $f_4$.
Then it is known that the minimal degree of a non-trivial Jacobian syzygy for $f_4 $ is $r-1$, see \cite[Example 4.5]{3syz}. It follows that one has either $r=2$ or $r=1$.

\medskip

\noindent {\bf Subcase 3A} $y$ occurs with a non-zero coefficient in $(a_1)_1$, say one has $(a_1)_1=\be x + \gamma y$, with $\gamma \ne 0$.
If  we write again $f$ as in \eqref{syz2}, we have
$$a_1f_x=\gamma p A_p(y)yx^{p-1} + \text{ higher order terms in } x,$$
$$b_1f_y=b_1A_p'(y)x^p+\text{ higher order terms in } x$$
and
$$xf=A_p(y)x^{p+1}+\text{ higher order terms in } x.$$
Hence in the sum 
\begin{equation}
\label{eqK4}
a_1f_x+b_1f_y+xf=0, 
\end{equation}
the term $\gamma p A_p(y)yx^{p-1}$
 cannot cancel out, a contradiction.
 
 We assume in the sequel in this Case that  $(a_1)_1=\be x$ with $\be \in \C$ and
 $(b_1)_1=\delta x + \gamma y$, with $\gamma \in \C^*$.

\noindent {\bf Subcase 3B} Assume in addition that $r=1$, that is
$f_4=\ell^4$ for some linear form $\ell$ in $R$.  Given our assumption, 
the first equality in \eqref{syz10} shows that $y^4$ cannot occur in $f_4$, while the second equality in \eqref{syz10} shows that $xy^3$ occurs in $f_4$. Since $f_4=(ux+vy)^4$, for some $u,v \in \C$, not both zero, this leads to a contradiction.

\noindent {\bf Subcase 3C} Assume now that $r=2$. There are two possibilities. Assume first that  $f_4=\ell_1^2\ell_2^2$, where $\ell_1$ and $\ell_2$ are linearly independent linear forms in $R$.
Then it follows from Arnold's classification lists, see \cite[Section 15.2]{AGV}, the case of singularities of corank 2 and zero 3-jet, $k=1$ on page 249, that $f$ is either homogeneous of degree 4, or a $T_{2,p,q}$ singularity. The homogeneous singularities are excluded form our discussion, while the $T_{2,p,q}$ singularity can be given by the equation
$$f=x^2y^2+x^p+y^q=0,$$
where $q \geq p \geq 4$ and $\frac{1}{p}+\frac{1}{q} <\frac{1}{2}$.
Such a singularity satisfies
$\mu(f)-\tau(f)=1$, and hence $f$ does not satisfy.
$$\mu(f)-\tau(f)=\dim_{\C} \frac{R}{ \Ann(f) }=2.$$
The remaining case to discuss is when $f_4=\ell_1^3\ell_2$, where $\ell_1$ and $\ell_2$ are linearly independent linear forms in $R$.
Note that
$$(f_4)_x=3\ell_1^2(\ell_1)_x\ell_2+\ell_1^3(\ell_2)_x=\ell_1^2(3(\ell_1)_x\ell_2+\ell_1(\ell_2)_x)$$
and
$$(f_4)_y=3\ell_1^2(\ell_1)_y\ell_2+\ell_1^3(\ell_2)_y=\ell_1^2(3(\ell_1)_y\ell_2+\ell_1(\ell_2)_y).$$
Hence the minimal Jacobian syzygy in this case is
$$(3(\ell_1)_y\ell_2+\ell_1(\ell_2)_y)(f_4)_x-(3(\ell_1)_x\ell_2+\ell_1(\ell_2)_x)(f_4)_y=0.$$
It follows that the pair $((a_1)_1,(b_1)_1)$ is proportional to the pair
$$(3(\ell_1)_y\ell_2+\ell_1(\ell_2)_y, -3(\ell_1)_x\ell_2-\ell_1(\ell_2)_x).$$
If we set
$$\ell_1=u_1x+v_1x \text{ and } \ell_2=u_2x+v_2y,$$
then a direct computation shows that we can assume that
\begin{equation}
\label{eqK5}
(a_1)_1=(u_1v_2+3u_2v_1)x+4v_1v_2y \text{ and } (b_1)_1=-4u_1u_2x-(3u_1v_2+u_2v_1)y,
\end{equation}
since the proportionality constant can be included in the coefficients of $\ell_1$ and $\ell_2$. By our assumption, $y$ does not occur in $(a_1)_1$, and hence $v_1v_2=0$.

\noindent {\bf Subcase 3Ca} 
Assume that $v_1=0$. Then $u_1 \ne 0$ and $v_2 \ne 0$ and hence
$$f_4=u_1^3x^3(u_2x+v_2y)=u_1^3u_2x^4+u_1^3v_2x^3y.$$
But  the second equality in \eqref{syz10} shows that $xy^3$ must occur in $f_4$, a contradiction.

\noindent {\bf Subcase 3Cb} 
Assume that $v_2=0$. Then $u_2 \ne 0$ and $v_1 \ne 0$ and hence
$$f_4=(u_1x+v_1y)^3u_2y=v_1^3u_2y^4+ \ldots$$
But the first equality in \eqref{syz10} shows that $y^4$ cannot occur in $f_4$, a contradiction.

\medskip

\noindent {\bf Case 4} $\ord f >4$.
Then \eqref{syz10} implies that $\ord a_2 \geq 2$ and $\ord b_2 \geq 2$. It also implies that $y$ cannot occur with a non-zero coefficient in either $(a_1)_1$ or in $(b_1)_1$ since otherwise one of the derivatives $f_x$ or $f_y$ would have order 3, a contradiction with $\ord f >4$. This ends the proof of Theorem \ref{thm3}.

\begin{ex}
\label{ex3}

(i) If we take
$$f= x^3+xy^8+y^{13},$$
then $\ord f=3$ and we get by a direct computation with SINGULAR \cite{Sing}
$$
\mu(f) = 22, \  \tau(f) = 20 \text{ and } \Ann(f)=\Fitt_2(I_f)=(x,y^2).$$

(ii) If we take
$$f= x^3y+xy^6+y^8,$$
then $\ord f=4$, $f_4=x^3y$ and we get by a direct computation 
$$
\mu(f) = 17, \  \tau(f) = 15 \text{ and } \Ann(f)=\Fitt_2(I_f)=(x,y^2).$$

Similarly, if we take
$$f= x^4+x^3y^2+x^3y^3+x^2y^4+xy^6,$$
then $\ord f=4$, $f_4=x^4$ and we get by a direct computation 
$$
\mu(f) = 21, \  \tau(f) = 19 \text{ and } \Ann(f)=\Fitt_2(I_f)=(x,y^2).$$

\end{ex}

The following example shows that the case $\dim_{\C} \frac{R}{ \Ann(f) }=3$ might be much more subtle.

\begin{ex}
\label{ex4}

 If we take
$$f= x^6+x^5y+y^7,$$
 we get by a direct computation with SINGULAR \cite{Sing}
$$
\mu(f) = 29, \  \tau(f) = 26 \text{ and } \Ann(f)=(6x+5y,y^3) \subset \Fitt_2(I_f)=(6x+5y,y^2).$$
On the other hand, if we take
$$f= x^3+xy^9+y^{14},$$
then $\ord f=3$ and we get by a direct computation 
$$
\mu(f) = 25, \  \tau(f) = 22 \text{ and } \Ann(f)=\Fitt_2(I_f)=(x,y^3).$$

\end{ex}



\end{document}